\documentclass[10pt,twoside]{article}
\usepackage{amsfonts}
\usepackage{amssymb}
\input xypic
\xyoption{all}
\LaTeXdiagrams
\xyoption{2cell}
\UseAllTwocells

\newtheorem{teo}{Theorem}

\newtheorem{rem}{Remark}
\def\be{\begin{equation}}
\def\ee{\end{equation}}
\def\lra{\longrightarrow}

\newenvironment{apl}[1]{\vspace{.5cm}{\bf \noindent Application: #1.}}{$\Box$}
\newenvironment{ins}[1]{\vspace{.5cm}{\bf \noindent Instance: #1.}}{$\Box$}
\newcounter{examnum}[section]
\newcounter{remarnum}[section]
\newcommand{\sqr}[6]{$$\xymatrix{ #1 \times #1
\ar[rr]^#3 && #2 \ar[dd]^#4
\\
\\
#1\ar[uu]^#6\ar[rr]_#5&&#2. }$$}
\begin{document}
\title{A Universal Approach to Self-Referential Paradoxes, Incompleteness and Fixed Points}
\author{Noson S. Yanofsky}
\date{}
\maketitle
\begin{flushright}
{\it The point of these observations\\
is not the reduction of the \\ familiar to the unfamiliar[...]\\
but the extension of the familiar\\ to cover many more cases.}\\
{Saunders MacLane}\\
{\it Categories for the Working Mathematician
}\cite{MR2001j:18001}\\ Page 226.
\end{flushright}
\begin{abstract}

\noindent Following F. William Lawvere, we show that many
self-referential paradoxes, incompleteness theorems and fixed
point theorems fall out of the same simple scheme. We demonstrate
these similarities by showing how this simple scheme encompasses
the semantic paradoxes, and how they arise as diagonal arguments
and fixed point theorems in logic, computability theory,
complexity theory and formal language theory.
\end{abstract}
\section{Introduction}

In 1969, F. William Lawvere wrote a paper \cite{MR39:4075} in
which he showed how to describe many of the classical paradoxes
and incompleteness theorems in a categorical fashion. He used the
language of category theory (and of cartesian closed categories in
particular) to describe the setting. In that paper he showed that
in a cartesian closed category satisfying certain conditions,
paradoxical phenomena can occur. Lawvere then went on to
demonstrate this scheme by showing the following examples
\begin{enumerate}
\item Cantor's theorem that $\mathbb{N} \lneqq \wp(\mathbb{N})$
\item Russell's paradox
\item The non-definability of satisfiability
\item Tarski's non-definability of truth and
\item G\"{o}del's first incompleteness theorem.
\end{enumerate}

Further work along these lines were done in several papers e.g.
\cite{MR91g:18007,MR91c:03049,MR93i:03021,MR91k:03018}.
Unfortunately, Lawvere's paper has been overlooked by many people
both inside and outside of the category theory community. Lawvere
and Schanuel revisited these ideas in Session 29 of their book
\cite{MR93m:18001}. Recently, Lawvere and Robert Rosebrugh came
out with a book {\it Sets for Mathematics} \cite{setsformath}
which also has a few pages on this scheme.

It is our goal to make these amazing results available to a larger
audience. Towards this aim we restate Lawvere's theorems without
using the language of category theory. Instead, we use sets and
functions. The main theorems and their proofs are done at tutorial
speed. We generalize one of the theorems and then we go on to show
different instances of these result. In order to demonstrate the
ubiquity of the theorems, we have tried to bring examples from
many diverse areas of logic and theoretical computer science.

Classically, Cantor proved that there is no onto (surjection)
function
$$\mathbb{N} \lra \bf 2^\mathbb{N}=\wp(\mathbb{N})$$ where
$\bf 2^\mathbb{N}$ is the set of functions from $\mathbb{N}$ to
${\bf 2} = \{0,1\}$. $2^\mathbb{N}$ is the set of characteristic
functions on the set $\mathbb{N}$ and is equivalent to the
powerset of $\mathbb{N}$. We can generalize Cantor's theorem to
show that for any set $T$ there is no onto function
$$ T \lra {\bf 2}^T = \wp(T).$$ The same theorem is also true for
other sets besides $\bf 2$, e.g. ${\bf 3}=\{0,1,2\}$ or ${\bf
23}=\{0,1,2, \ldots 21,22\}$. The theorem is not true for the set
${\bf 1}=\{0\}$. In general we can replace $\bf 2$ with an
arbitrary ``non-degenerate'' set $Y$. From this generalization,
the basic statement of Cantor's theorem roughly says that if $Y$
is ``non-degenerate'' then there is no onto function
$$T \lra Y^T$$ where $Y^T$ is the set of functions from $T$ to
$Y$. $Y$ can be thought of as the set of possible ``truth-values''
or ``properties'' of elements of $T$. By ``non-degenerate'' we
mean that the objects of $Y$ can be interchanged or that there
exists a function $\alpha$ from $Y$ to $Y$ without any fixed
points ($y\in Y$ where $\alpha(y)=y$.)

 Rather than looking at functions $\widehat{f}:T \lra Y^T$,
 we shall look at equivalent functions of
the form $f:T\times T \lra Y$. Every $\widehat{f}$ can be
converted to a function $f$ where $f(t,t')=\widehat{f}(t')(t) \in
Y$. Saying that $\widehat{f}$ is not onto is the same thing as
saying that there exists a $g(-) \in Y^T$ such that for all $t'
\in T$ the function $\widehat{f}(t')=f(-,t'):T \lra Y$ is not the
same as the function $g(-):T \lra Y$. In other words there exists
a $t \in T$ such that
$$g(t) \neq f(t,t').$$ We shall call a function $g:T \lra Y$
``representable by $t_0$'' if $g(-)=f(-,t_0)$. So if $\widehat{f}$
is not onto, then there exists a $g(-) \in Y^T$ that is not
representable by any $t\in T$.

On a philosophical level, this generalized Cantor's theorem says
that as long as the truth-values or properties of $T$ are
non-trivial, there is no way that a set $T$ of things can ``talk
about'' or ``describe'' their own truthfulness or their own
properties. In other words, there must be a limitation in the way
that $T$ deals with its own properties. The Liar paradox is the
three thousand year-old primary example that shows that natural
languages should not talk about their own truthfulness. Russell's
paradox shows that naive set theory is inherently flawed because
sets can talk about their own properties (membership.) G\"{o}del's
incompleteness results shows that arithmetic can not talk
completely about its own provability. Turing's Halting problem
shows that computers can not completely deal with the property of
whether a computer will halt or go into an infinite loop. All
these different examples are really saying the same thing: there
will be trouble when things deal with their own properties. It is
with this in mind that we try to make a single formalism that
describes all these diverse -- yet similar -- ideas.

The best part of this unified scheme is that it shows that there
are really no paradoxes. There are limitations. Paradoxes are ways
of showing that if you permit one to violate a limitation, then
you will get an inconsistent systems. The Liar paradox shows that
if you permit natural language to talk about its own truthfulness
(as it - of course - does) then we will have inconsistencies in
natural languages. Russell's paradox shows that if we permit one
to talk about any set without limitations, we will get an
inconsistency in set theory.  This is exactly what is said by
Tarski's theorem about truth in formal systems. Our scheme shows
the inherent limitations of all these systems. The constructed
$g$, in some sense is the limitation that your system ($f$) can
not deal with. If the system does deal with the $g$, there will be
an inconsistency (fixed point).

The contrapositive of Cantor's theorem says that if there is a
onto $T \lra Y^T$ then $Y$ must be ``degenerate'' i.e. every map
from $Y$ to $Y$ must have a fixed point. In other words, if $T$
can talk about or describe its own properties then $Y$ must be
faulty in some sense. This ``degenerate''-ness is a way of
producing fixed point theorems.

For pedagogical reasons, we have elected not to use the powerful
language of category theory. This might be an error. Without using
category theory we might be skipping over an important step or
even worse: wave our hands at a potential error. It is our hope
that this paper will make you go out and look at Lawvere's
original paper and his subsequent books. Only the language of
category theory can give an exact formulation of the theory and
truly encompass all the diverse areas that are discussed in this
paper. Although we have chosen not to employ category theory here,
its spirit is nevertheless pervasive throughout.

This paper is intended to be extremely easy to read. We have tried
to make use of the same proof pattern over and over again.
Whenever possible we use the same notation. The examples are
mostly disjoint. If the reader is unfamiliar with or can not
follow one of them, he or she can move on to the next one without
losing anything. Section 2 states Lawvere's main theorem and some
of our generalizations. Section 3 has many worked out examples. We
start the section with the classical paradoxes and then move on
some of the semantic paradoxes. From there we go on to other
examples from theoretical computer science. Section 4 states the
contrapositive of the main theorem and some of its
generalizations. The examples of this contrapositives are in
Section 5. We finish off the paper by looking at some future
directions for this work to continue. We also list some other
examples of limitations and fixed point theorems that might be
expressible in our scheme.

We close this introduction with a translation of Cantor's original
proof of his diagonalization theorem. His language is remarkably
reminiscent of our language. This translation was taken from
Shaughan Lavine's book \cite{MR95k:00009}.
\begin{quotation}
The proof seems remarkable not only because of its simplicity, but
especially also because the principle that is employed in it can
be extended to the general theorem, that the powers of
well-defined sets have no maximum or, what is the same, that for
any given set $L$ another $M$ can be placed beside it that is of
greater power than $L$.

For example Let $L$ be a linear continuum, perhaps the domain of
all real numerical quantities that are $\geq 0$ and $\leq 1$.

Let $M$ be understood as the domain of all single-valued functions
$f(x)$ that take on only the two values $0$ or $1$, while $x$ runs
through all real values that are $\geq 0$ and $\leq 1$. [
$M=\mathbf{2}^L$...]

But $M$ does not have the same power as $L$ either. For otherwise
$M$ can be put into one-to-one correspondence to the variable $z$
[of $L$], and thus $M$ could be thought of in the form of a single
valued function
$$\phi(x,z)$$
of the two variables $x$ and $z$, in such a way that through every
specification of $z$ one would obtain an element $f(x)=\phi(x,z)$
of $M$ and also conversely each element $f(x)$ of $M$ could be
generated from $\phi(x,z)$ through a single definite specification
of $z$. This however leads to a contradiction. For if we
understand by $g(x)$ that single valued function of $x$ which
takes only values $0$ or $1$ and which every value of $x$ is
different from $\phi(x,x)$, then on the one hand $g(x)$ is an
element of $M$, and on the other it can not be generated from
$\phi(x,z)$ by any specification $z=z_0$, because $\phi(z_0,z_0)$
is different from $g(z_0)$.
\end{quotation}

{\bf Acknowledgments.} The  author is grateful to Rohit Parikh for
suggesting that this paper be written and for his warm
encouragement. The author also had many helpful conversations with
Eva Cogan, Scott Dexter, Mel Fitting, Alex Heller, Roman Kossak,
Mirco Mannucci, and Paula Whitlock.

\section{Cantor's Theorems and its Generalizations}
It is pedagogically sound to skip this section for a moment and
read the beginning of the next section where you can remind
yourself of the proof of the more familiar version of Cantor's
theorem (about $\mathbb{N} \lneqq \wp(\mathbb{N})$) and Russell's
set theory paradox. Our theorem here might seem slightly abstract
at first.

\begin{teo}[Cantor's Theorem] If $Y$ is a set and there exists a function $\alpha:Y\lra Y$
without a fixed point (for all $y \in Y$, $\alpha(y)\neq y$), then
for all sets $T$ and for all functions $f:T \times T \lra Y$ there
exists a function $g:T \lra Y$ that is not representable by $f$
i.e.  such that for all $t \in T$ $$g(-) \neq f(-, t).$$
\end{teo}

\noindent{\bf Proof.} Let $Y$ be a set and assume $\alpha:Y \lra
Y$ is a function without fixed points. There is a function
$\bigtriangleup:T \lra T\times T$ that sends every $t \in T$ to
$(t,t) \in T\times T$. Then construct $g:T \lra Y$ as the
following composition of three functions.
\sqr{T}{Y}{f}{\alpha}{g}{\bigtriangleup} In other words, $$ g(t) =
\alpha(f(t,t)).$$ We claim that for all $t \in T$, $g(-) \neq f(-,
t)$ as functions of one variable. If $g(-) = f(-, t_0)$ then by
evaluation at $t_0$ we have
$$f(t_0,t_0)= g(t_0) = \alpha(f(t_0,t_0))$$
where the first equality is the fact that $g$ is representable and
the second equality is the definition of $g$. But this means that
$\alpha$ does have a fixed point. $\Box$

\begin{rem} Obviously, every set with two or more elements has a
function to itself that does not have a fixed point. It is here
that we get in trouble for talking about sets and functions as
opposed to objects in a category and morphisms between those
objects. Perhaps $Y$ and $T$ are sets with extra (algebraic)
structure and functions between them are intended to preserve that
extra structure. In that case, we are really dealing with fewer
functions between the sets.
\end{rem}

\begin{rem} The $\bigtriangleup$ map is called the ``diagonal''
and many of the proofs are called ``diagonalization arguments.''
$f$ is some type of evaluation function and $f(t,t)$ is an
evaluation of itself, hence ``self-reference'' or
``self-referential arguments.''
\end{rem}

\begin{rem} We follow Lawvere and Schanuel \cite{MR93m:18001} in calling this theorem ``Cantor's
Theorem'' and it's contrapositive the ``Diagonal Theorem'' stated
in Section 4.
\end{rem}

We generalize the above theorem so that instead of $\bigtriangleup
= \langle Id,Id \rangle$ we use $\langle Id,\beta\rangle$ for an
arbitrary onto (right invertible) function $\beta:T \lra S$.
Whereas $\bigtriangleup = \langle Id,Id \rangle :T \lra T \times
T$ takes every $t$ to $(t,t)$, $\langle Id,\beta\rangle :T \lra T
\times S$ takes every $t$ to $(t,\beta(t))$.

The way to think about this theorem is to say that if there is a
onto $\beta: T \lra S$ then in a sense $|S| \leqslant |T|$ and
Cantor's theorem says $|T| \leqslant |Y^T|$ and so we conclude
that $|S| \leqslant |Y^T|$.

\begin{teo} Let $Y$ be a set, $\alpha:Y\lra Y$ a function
without a fixed point, $T$ and $S$ sets and $\beta:T \lra S$ a
function that is onto (i.e., has a right inverse $\bar{\beta}:S
\lra T$,) then for all functions $f:T \times S \lra Y$ the
function $g_\beta:T \lra Y$ constructed as follows
$$\xymatrix{ T \times S
\ar[rr]^f && Y \ar[dd]^\alpha
\\
\\
T\ar[uu]^{\langle Id,\beta\rangle}\ar[rr]_{g_\beta}&&Y. }$$
 is not
representable by $f$.
\end{teo}
{\bf Proof.} Let $Y,\alpha, T$ and $\beta$ be given. Let
$\bar{\beta}:S \lra T$ be the right inverse of $\beta$. By
definition
$$g_\beta(t) = \alpha(f(t,\beta(t))).$$
We claim that for all $s \in S$ $g_\beta(-) \neq f(-, s).$ If
$g_\beta(-) = f(-, s_0)$ then evaluation at $\bar{\beta}(s_0)$
gives
\begin{eqnarray*}
f(\bar{\beta}(s_0),s_0)& = & g_\beta(\bar{\beta}(s_0))
\quad \mbox{ by representability of } g_\beta\\
& = & \alpha(f(\bar{\beta}(s_0),\beta(\bar{\beta}(t_0)))) \quad \mbox{ by definition of } g_\beta\\
& = & \alpha(f(\bar{\beta}(s_0),s_0)) \quad \mbox{ by definition of right inverse.} \\
\end{eqnarray*} Which means that $\alpha$ does have a fixed point.
 $\Box$
\newline \newline

We can think of this theorem in another way. Set $S=T$ and lets
consider a $\beta$ different than $Id_T$. The usual way to
visualize Cantor's Theorem is
$$ \begin{array}{l|lllllc}
f&t_1&t_2&t_3&t_4&t_5& \cdots \\
\hline
t_1&[y_3]&y_7& y_{21}&y_2&y_4&\cdots\\
t_2&y_1&[y_{17}]& y_{2}&y_7&y_{41}&\cdots \\
t_3&y_0&y_3& [y_{7}]&y_2&y_{24}&\cdots \\
t_4&y_9&y_7& y_{64}&[y_2]&y_4&\cdots \\
t_5&y_4&y_{73}& y_{31}&y_2&[y_4]&\cdots \\
\vdots &\vdots& & \vdots& &\vdots &\ddots \\
\end{array}
$$
Everything that is in square brackets gets changed. For example
$y_3$ gets changed to $\alpha(y_3)$.  However a little thought
shows that we do not need to go along the diagonal. The diagonal
is just the simplest way. What is needed is that every row of the
table gets at least one element changed. So we might have a
picture that looks like this:
$$ \begin{array}{l|lllllc}
f&t_1&t_2&t_3&t_4&t_5& \cdots \\
\hline
t_1&y_3&y_7& y_{21}&[y_2]&y_4&\cdots\\
t_2&[y_1]&y_{17}& [y_{2}]&y_7&y_{41}&\cdots \\
t_3&y_0&y_3& y_{7}&y_2&[y_{24}]&\cdots \\
t_4&y_9&[y_7]& y_{64}&[y_2]&y_4&\cdots \\
t_5&y_4&y_{73}& y_{31}&[y_2]&y_4&\cdots \\
\vdots &\vdots& & \vdots& &\vdots &\ddots \\
\end{array}
$$
The fact that every row has something changed is in essence the
fact that $\beta$ is onto. As long as $\beta$ is onto, Cantor's
theorem still holds.

With this in mind we may pose -- but do not answer -- the
following questions. Should these theorems really be called
``diagonalization theorems''? Does self-reference really play a
role here? Since we can generate the same paradoxes without
self-reference, does this destroy Russell's vicious-circle
principle?

\section{Instances of Cantor's Theorems}
We shall begin with the familiar version of Cantor's theorem about
the power set of the natural numbers. From there we move on to
Russell's set theory paradox and other paradoxes and limitations.
We shall do the first two instances slowly and use the same
notation and ideas as the theorems in the last section. The other
instances we shall do more quickly.

\begin{ins}{Cantor's $\mathbb{N} \lneqq \wp(\mathbb{N})$ Theorem} The
theorem says that there can not be an onto function from
$\mathbb{N}$ to $\wp(\mathbb{N})$. Let $S_0, S_1, S_2, \ldots$ be
a proposed enumeration of all subsets of $\mathbb{N}$. Let ${\bf
2}=\{0,1\}$ be a set and consider the ``negation'' function
$\alpha:{\bf 2} \lra {\bf 2}$ where $\alpha(0)=1$ and
$\alpha(1)=0$. Let $f:\mathbb{N} \times \mathbb{N} \lra \bf 2$ be
defined as

$$ f(n,m)= \left \{
\begin{array}{r@{\quad : \quad}l}
1 & \mbox{if } n \in S_m\\
0 & \mbox{if } n \not\in S_m. \end{array}\right.$$ For each $m$,
$f(-,m)$ is the characteristic function of $S_m$: $$f(-,m) =
\chi_{S_m}.$$ Construct $g$ as follows: \sqr{\mathbb{N}}{{\bf
2}}{f}{\alpha}{g}{\bigtriangleup} $g$ is the characteristic
function of the set $$G=\{n \in \mathbb{N} | n \not\in S_n \}.$$
For all $m$, $\chi_G=g(-) \neq f(-, m)=\chi_{S_m}$. Because if
there was an $m_0$ such that $g(-)= f(-, m_0)$ then by evaluation
at $m_0$ we have
$$f(m_0,m_0)= g(m_0) = \alpha(f(m_0,m_0))$$
where the first equality is from the fact that $g$ is
representable by $m_0$ and the second equality is by the
definition of $g$. This means that the negation operator has a
fixed point which is clearly false. In other words $G \subseteq
\mathbb{N}$ is not in the proposed enumeration of all subsets of
$\mathbb{N}$.
\end{ins}

\begin{ins}{Russell's Paradox} This paradox says that the set of all
sets that are not members of themselves is both a member of itself
and not a member of itself. Let $Sets$ be some universe of sets
(we are being deliberately ambiguous here.) Again consider the
``negation'' function $\alpha:{\bf 2}\lra {\bf 2}$ where
$\alpha(0)=1$ and $\alpha(1)=0$. Let $f:Sets \times Sets \lra \bf
2$ be defined as follows on sets $s$ and $t$.

 $$ f(s,t)= \left \{
\begin{array}{r@{\quad : \quad}l}
1 & \mbox{if } s \in t\\
0 & \mbox{if } s \not\in t. \end{array}\right.$$  We construct
$g$ as follows
 \sqr{Sets}{{\bf 2}}{f}{\alpha}{g}{\bigtriangleup}
$g$ is the characteristic function of those sets that are not a
member of themselves. For all sets $t$, $g(-) \neq f(-, t)$.
Because if there was a set $t_0$ such that $g(-)= f(-, t_0)$ then
from evaluation at $t_0$ we get
$$f(t_0,t_0)= g(t_0) = \alpha(f(t_0,t_0))$$
where the first equality is because $g$ is representable and the
second equality is from the definition of $g$. This is plainly
false. To summarize, in order to make sure that there are no
paradoxes we must say that $g$ is the characteristic function of a
``collection'' of $Sets$ but this ``collection'' does not form a
set.

We mention in passing that the Barber paradox and other simple
self-referential paradoxes can be done exactly like this. The
Barber paradox has a simple solution, namely that the village
described by the phrase ``there is a village where everyone who
does not shave themselves is shaved by the barber'' does not
really exist. We are in a sense saying the same thing about
Russell's paradox. Namely, the collection of sets that do not
contain themselves does not form an existent set. For some reason,
people find it more ontologically disheartening to say that a
collection does not form a set than that a particular village does
not exist.
\end{ins}

\begin{ins}{Grelling's Paradox} We now move on to some of the semantic paradoxes.
There are some adjectives that describe themselves and there are
some that do not. ``English'' is an English word. ''French'' is
not a French word. ``Short'' is not short and ``Long'' is not
long. ``Polysyllabic'' is polysyllabic but ``monosyllabic'' is not
monosyllabic. Call all words that do not describe themselves
``heterological.'' Now ask yourself if ``heterological'' is
heterological. It is if and only if it is not.

Consider the set $Adj$ of all (English) adjectives. We have the
following function $f:Adj \times Adj \lra \bf 2$ defined for all
adjectives $a_1$ and $a_2$,
 $$ f(a_1,a_2)= \left \{
\begin{array}{r@{\quad : \quad}l}
1 & \mbox{if } a_2 \mbox{ describes } a_1 \\
0 & \mbox{if } a_2 \mbox{ does not decribe } a_1.
\end{array}\right.$$ And so we have the following construction of $g$

\sqr{Adj}{{\bf 2}}{f}{\alpha}{g}{\bigtriangleup} $g$ is the
characteristic function of a subset ( = property) of adjectives
that can not be described by any adjectives. This is exactly what
is meant by $g(-) \neq f(-,a)$ for all adjectives $a$.
``Heterological'' is not the only adjective that is in this
subset. Some authors (e.g. Kleene) have also used the word
``impredicable''. Our formulation includes all such paradoxical
adjectives.
\end{ins}

\begin{ins}{Liar Paradox}  The oldest example of a self-referential
paradox is the (Cretans) liar paradox.  Epimenides of Crete said
``All Cretans are liars.'' There are many such examples: ``This
sentence is false.'', ``I am lying.'' The Liar paradox is very
similar to Grelling's paradox. Whereas with Grelling's paradox we
dealt with adjectives, here we deal with complete English
sentences. Quine's paradox is the primary example:
\begin{quote} `yields
falsehood when appended to its own quotation'\\
yields falsehood when appended to its own quotation.\end{quote}
The philosophical literature is full of such examples. Since the
formalism is similar to Grelling's paradox, we leave it to the
reader.
\end{ins}

\begin{ins}{The Strong Liar Paradox} A common ``solution'' to the
Liar's paradox is to say that that there are certain sentences
that are neither true nor false but are meaningless. ``I am
lying'' would be such a sentence. This is a type of three-valued
logic. This is, however, not a ``solution.'' Consider the sentence
\begin{quote} `yields
falsehood or meaninglessness \\ when appended to its own quotation'\\
yields falsehood or meaninglessness \\ when appended to its own
quotation.\end{quote} If this sentence is true, then it is false
or meaningless. If it is false, then it is true and not
meaningless. If it is meaningless, then it is true and not
meaningless.

This paradox can also be formulated with our scheme. Consider the
set of English sentences $Sent$ and the set ${\bf 3} = \{T(rue),
M(eaningless), F(alse) \}$. We have the following function $f:Sent
\times Sent \lra \bf 3$ defined for all sentences $s_1$ and $s_2$,
 $$ f(s_1,s_2)= \left \{
\begin{array}{r@{\quad : \quad}l}
T & \mbox{if } a_2 \mbox{ describes } a_1 \\
M & \mbox{if } \mbox{ it is meaningless for } a_2 \mbox{ to describe } a_1 \\
F & \mbox{if } a_2 \mbox{ does not decribe } a_1.
\end{array}\right.$$ Now consider the function $\alpha: \bf 3 \lra \bf 3$ defined as
$\alpha(T)=F$ and $\alpha(M)=\alpha(F)=T$. Construct $g$ as
follows \sqr{Sent}{{\bf 3}}{f}{\alpha}{g}{\bigtriangleup} $g$ is
the characteristic function of sentences that are neither false
nor meaningless when describing themselves. By characteristic
function we mean those sentences that $g$ takes to $T$ as opposed
to $M$ or $F$.
\end{ins}

\begin{ins}{Richard's Paradox}
There are many sentences in the English language that describe
real numbers between $0$ and $1$.   Let us lexicographically order
all English sentences. Using this order, we can select all those
English sentences that describe real numbers between $0$ and $1$.
For example ``x is the ratio between the circumference and the
diameter of a circle divided by ten.'' describes the number
$0.314159 \ldots$. There are many similar English sentences. Call
such a sentence a ``Richard Sentence.'' So we have the concept of
the ``$m$-th Richard Sentence.''

Consider the set ${\bf 10}=\{0,1,2,\ldots 9\}$ and the function
$\alpha:{\bf 10} \lra {\bf 10}$ defined as $\alpha(i)= 9-i$. This
function does not have a fixed point. Now consider the function
$f:\mathbb{N} \times \mathbb{N} \lra {\bf 10}$ defined as
$$f(n,m)=\textrm{The $n$-th decimal number of the $m$-th Richard
Sentence.}$$ For example, if the sentence in the above paragraph
is the 15th Richard sentence then $f(4,15)=1$ because of the $1$
in $0.314{\bf 1}59 \ldots$. Now consider $g:\mathbb{N} \lra {\bf
10}$ constructed as  \sqr{\mathbb{N}}{{\bf
10}}{f}{\alpha}{g}{\bigtriangleup} This $g$ describes a real
number between $0$ and $1$ and yet for all $m \in \mathbb{N}$
$$ g(-) \neq f(-,m)$$ i.e. this number is different than all
Richard Sentences. Yet here is a Richard Sentence that describes
this number:
\begin{quote}
x is the real number between $0$ and $1$ whose $n$-th digit is
nine minus the $n$-th digit of the number described by the $n$-th
Richard sentence.
\end{quote}
For reasons that are beyond the author, this paradox remains.
\end{ins}

\begin{ins}{Turing's Halting Problem}  The following formulation was inspired
by Heller's fascinating work on recursion categories \cite{Heller}
and Manin's intriguing paper on classical and quantum computations
\cite{Manin}.

For this instance we leave the comfortable world of sets and
functions. We must talk about computable universes. A computable
universe is a category $\bf U$ with the following two properties

\begin{enumerate}
\item $\mathbb{N}$ and $\bf 2$ are objects in $\bf U$
\item For every object $C$ in $\bf U$ there is some type of
enumeration of the elements of $C$. An enumeration is a total
isomorphism $e_C:\mathbb{N} \lra C$. One should think of $C$ as a
set of computable things, e.g., trees, graphs, numbers, stacks,
strings etc.
\item For every (not necessarily total) function $f:C \lra C'$
there is a corresponding number $\ulcorner f \urcorner \in
\mathbb{N}$. Think of this as the G\"{o}del number of the program
that computes the computation.
\item For every (not necessarily
total) function $f:C \lra C'$ there is a corresponding recursively
enumerable (r.e.) set $W_{\langle f \rangle} \subseteq
\mathbb{N}$. For every $c \in C$, $f$ has a value at $c$ if and
only if $e_C^{-1}(c) \in W_{\langle f \rangle}$. Again one should
think of a partial function from one computable domain to another.
\end{enumerate}

$Halt$ in a computable universe should be a total function
$Halt:\mathbb{N} \times \mathbb{N} \lra \bf 2$ in $\bf U$ such
that for all $f:C\lra C'$
$$Halt(-,\ulcorner f \urcorner) = \chi_{W_{\langle f \rangle}}.$$
This says that $Halt$ should be able to tell for what values in
$C$ the computation halts. Formally
$$
Halt(n,m)= \left \{
\begin{array}{r@{\quad : \quad}l}
1 & \mbox{if } n \in W_m\\
0 & \mbox{if } n \not\in W_m. \end{array}\right.$$

Consider $\alpha: \bf 2 \lra \bf 2$ defined as follows:
$\alpha(0)=1$ and $\alpha(1)\uparrow$, i.e., the computation is
undefined. Construct $g$ as follows: \sqr{\mathbb{N}}{{\bf
2}}{{Halt}}{\beta}{g}{\bigtriangleup}

We conclude by showing that $Halt$ is not total because it is not
defined at $\ulcorner g \urcorner$. If $Halt$ was defined at
$\ulcorner g \urcorner$ then we would have the following
contradiction:
\begin{eqnarray*}
Halt(\ulcorner g\urcorner,\ulcorner g\urcorner) = 1 & \mbox{iff} &
\ulcorner g\urcorner \in W_{\langle g \rangle}
\quad \mbox{ by definition of } Halt \\
& \mbox{iff} & g(\ulcorner g\urcorner) = 1 \quad \mbox{ by the halting of } g\\
& \mbox{iff} & Halt(\ulcorner g\urcorner,\ulcorner g\urcorner) = 0
\quad \mbox{ by the definition of } g.\\
\end{eqnarray*}
Hence no total $Halt$ can exist.
\end{ins}

\begin{ins}{A non-r.e. Language} There is a language that is not
recognized by any Turing machine. Let $M_0, M_1 , M_2, \ldots$ be
an enumeration of all Turing machines on the input language
$\Sigma=\{0,1\}$. Let $w_0,w_1, w_2, \ldots $ be an enumeration of
all the words in $\Sigma^*$. If $w_i$ is a word in $\Sigma$ we let
$(w_i)$ denote the numerical value of the binary word. Consider
the following function $f:\Sigma^* \times \Sigma^* \lra \bf 2$
defined as follows: $$ f(w_i,w_j)= \left \{
\begin{array}{r@{\quad : \quad}l}
1 & \mbox{if } w_i \mbox{ is accepted by } M_{(w_j)} \\
0 & \mbox{if } w_i \mbox{ is not accepted by  } M_{(w_j)}.
\end{array}\right.$$ Then the constructed $g$
\sqr{\Sigma^*}{{\bf 2}}{f}{\alpha}{g}{\bigtriangleup} is the
characteristic function of a language that is not accepted by any
Turing machine. Of course, the fact that there are non-r.e.
languages also follows from a simple counting argument. Namely the
number of Turing machines is countable and the number of languages
($\wp(\Sigma^*)$) is uncountable.
\end{ins}

\begin{ins}{An Oracle $B$ such that $P^B \neq NP^B$} One of the
major open questions in computer science is whether or not $P$,
the set of all problems that can be solved by deterministic Turing
machines (TMs) in polynomial time, is equal to the set $NP$, of
all problems that can be solved by non-deterministic TMs in
polynomial time. Alas, this question will not be answered in this
paper. However there is a related question that can be answered.
Consider the same question for oracle TMs. An oracle TM is a TM
with an associated set $S$, such that the TM can determine if a
word is actually an element of $S$. For a given set $S$ there are
analogous sets $P^S$ and $NP^S$. Baker, Gil and Solovay
\cite{MR52:16108} have proven that there exists a set $A$ such
that $P^A = NP^A$ and there exists a set $B$ such that $P^B \neq
NP^B$. Here we shall prove the second result. Since every
deterministic machine is by definition also nondeterministic, we
have for every $B$,  $P^B \subseteq NP^B$. What remains is to show
that there is a set $B$ and a language $L_B$ such that $L_B \in
NP^B$ but $L_B \not\in P^B$ i.e. $NP^B \nsubseteq P^B$. Our proof
was adopted from \cite{MR83j:68002}.

Let $M_0^?, M_1^?, M_2^?, \ldots$ be some enumeration of all the
oracle deterministic polynomial Turing machines in the alphabet
$\Sigma = \{0,1\}.$ There is a corresponding sequence of
polynomials $p_0(x), p_1(x), p_2(x), \ldots$ expressing the worst
execution time for each machine.

For any function $f:\Sigma^* \times \mathbb{N} \lra \bf 2$ and for
each $i\in \mathbb{N}$, $f(-, i):\Sigma^* \lra \bf 2$ is a
characteristic function on the set $\Sigma^*$. We will often
confuse a set and its characteristic function. Let
$\overline{f}(-,i)$ denote the characteristic function of the
complement of $f(-,i)$, i.e., $\overline{f}(-,i)$ is the set that
$f(-,i)$ takes to $0$. Let $\overline{F}(-,i)$ denote the
cumulative characteristic function
$$\overline{F}(-,i) = \bigcup_{j\leq i} \overline{f}(-,j).$$

We shall define $f(-,-)$ inductively. $(\forall w \in \Sigma^*)
f(w,0)=1.$ For $w\in \Sigma^*$ and $i \in \mathbb{N}$, $f(w,i)=0$
if and only if the following three conditions are satisfied
\begin{enumerate}
\item $(\forall w'<w)  f(w',i)=1$ where the $<$ is a lexicographical
order on the words of $\Sigma^*$. This insures that there is only
one word accepted to $B$ for each $i$.
\item $M_i^{\overline{F(-,i)}}$ rejects $0^{|w|}$ within  $i^{log\
i}$ steps.
\item $(\forall j<i) M_j^{\overline{F(-,j)}}$ on input $0^{|w|}$
does not to query $w$ within  $j^{log\ j}$ steps.
\end{enumerate}
Once this $f$ is defined, we construct $g$ as follows
$$\xymatrix{ \Sigma^* \times \mathbb{N} \ar[rr]^f && \bf 2
\ar[dd]^\alpha
\\
\\
\Sigma^* \ar[uu]^{\langle Id,\beta \rangle }\ar[rr]_{g_\beta}&&
\bf 2 }$$ where $\beta(w)=|w|$, $\alpha(0)=1$ and $\alpha(1)=0$.
$g(w)=1$ if and only if $f(w,|w|)=0$ if and only if the above
three requirements are satisfied.

$g$ is the characteristic function of the set $B\subseteq
\Sigma^*$.  Now construct the language $$L_B= \{ 0^i | B \mbox{
contains a word of length } i \}.$$ This language can easily be
recognized by a linear time nondeterministic TM. On input $0^i$,
the NTM simply has to guess a string $w$ of length $i$ and see if
it is in $B$. Hence $L_B \in NP^B$. In contrast, because of
condition 2 above, $L_B$ can not be recognized by any DTM in
polynomial time, i.e., $(\forall m) g(-) \neq f(-,m).$
\end{ins}

\section{Diagonal Theorem and Generalizations}
The contrapositive of Cantor's Theorem is of equal importance.
\begin{teo}[Diagonal Theorem] If $Y$ is a set and there exists a set $T$
and   a function $f:T\times T \lra Y$ such that all functions $g:T
\lra Y$ are representable by $f$ (there exists a $t \in T$ such
that $g(-) = f(-,t),$) then all functions $\alpha:Y \lra Y$ have a
fixed point.
\end{teo} {\bf Proof.} The proof is constructive. Let $Y, T, f$
and $\alpha$ be given. Then we construct $g$ as follows:
\sqr{T}{Y}{f}{\alpha}{g}{\bigtriangleup}
 $g$ is defined as $$g(m) = \alpha(f(m,m)).$$ Since we have
assumed that $g$ is representable by some $t\in T$, we have that
$$g(m) = f(m,t).$$ And so we have a fixed point of $\alpha$ at
$y_0=g(t)$. Explicitly we have
\begin{eqnarray*}
\alpha(g(t)) & = & \alpha(f(t,t))
\quad \mbox{ by representation of } g\\
& = & g(t) \quad \mbox{ by definition of } g\\
\end{eqnarray*} $\Box$
\begin{rem}
Obviously, any set $Y$ with two or more elements has functions $Y
\lra Y$ that do not have fixed points. It is here that we get in
trouble by ignoring the category theory that is necessary. In the
examples that we will do, the objects we will be dealing with have
more structure then just sets and the functions between the
objects are required to preserve that structure. We are only
talking about these restricted functions.
\end{rem}

\begin{rem} It is important to note that the theorem uses a stronger hypothesis
than the proof actually uses. The theorem asks that {\bf all} $g:T
\lra Y$ be representable, however the proof only uses the fact
that any $g$ constructed in such a manner is representable. In the
future, we shall use this fact and only require that constructed
$g$ be representable.
\end{rem}

\section{Instances of Diagonal Theorems}
We use Mendelson's \cite{MR99b:03002} notation and language. In
particular $ \ulcorner\mathcal{B}(x)\urcorner$ is the G\"{o}del
number of $ \mathcal{B}(x)$. We shall assume that we are working
in a theory where there is a recursive $D:\mathbb{N} \lra
\mathbb{N}$ that is defined as follows: For all $\mathcal{B}(x)$
where $\mathcal{B}$ is a logical statement with $x$ its only free
variable then $$ D(\ulcorner\mathcal{B}(x)\urcorner) =
\ulcorner\mathcal{B}(\ulcorner\mathcal{B}(x)\urcorner)\urcorner.$$
\begin{teo}[Diagonalization Lemma]
For any well-formed formula (wf) $\mathcal{E}(x)$  with $x$ as its
only free variable, there exists a closed formula $\mathcal{C}$
such that $$\vdash \mathcal{C} \longleftrightarrow
\mathcal{E}(\ulcorner \mathcal{C} \urcorner). $$
\end{teo}
{\bf Proof.} Let $Lind^i$ be the set of Lindenbaum classes
(algebra) of well-formed formulas with $i$ free variables. Two wfs
are equivalent iff they are provably logically equivalent. Let
$f:Lind^1 \times Lind^1 \lra Lind^0$ be defined for two wfs with a
free variable $\mathcal{B}(x)$ and $\mathcal{H}(y)$ as follows:
$$f(\mathcal{B}(x), \mathcal{H}(y)) =
\mathcal{H}(\ulcorner\mathcal{B}(x)\urcorner).$$ Let the operator
on $Lind^0$ $\Phi_\mathcal{E}:Lind^0 \lra Lind^0$ be defined as
$\mathcal{P} \mapsto \Phi_\mathcal{E}(\mathcal{P}) =
\mathcal{E}(\ulcorner\mathcal{P}\urcorner)$. Using these
functions, we combine them to create $g$ as follows:
\sqr{Lind^1}{Lind^0}{f}{{\Phi_{\mathcal{E}}}}{g}{\bigtriangleup}
By definition $$g(\mathcal{B}(x))=\Phi_\mathcal{E}(
f(\mathcal{B}(x), \mathcal{B}(x))) =
\mathcal{E}(\ulcorner\mathcal{B}(\ulcorner\mathcal{B}(x)\urcorner)\urcorner).$$
We claim that $g$ is representable by $\mathcal{G}(x) =
\mathcal{E}(D(x))$. This is true because $$
g(\mathcal{B}(x))=\mathcal{E}(\ulcorner\mathcal{B}(\ulcorner\mathcal{B}(x)\urcorner)\urcorner)=
\mathcal{E}(D(\ulcorner\mathcal{B}(x)\urcorner))=\mathcal{G}(\ulcorner\mathcal{B}(x)\urcorner)=
f(\mathcal{B}(x), \mathcal{G}(y)). $$ So there is a fixed point of
$\Phi_\mathcal{E}$ at $\mathcal{C} =
\mathcal{G}(\ulcorner\mathcal{G}(x)\urcorner)$. Explicitly we have
\begin{eqnarray*}
\mathcal{E}(\ulcorner\mathcal{G}(\ulcorner\mathcal{G}(x)\urcorner)\urcorner)
& = &
\Phi_\mathcal{E}(\ulcorner\mathcal{G}(\ulcorner\mathcal{G}(x)\urcorner)\urcorner)
\quad \mbox{ by definition of } \Phi_\mathcal{E} \\
& = & \Phi_\mathcal{E}(f(\mathcal{G}(x),\mathcal{G}(x))) \quad \mbox{ by definition of } f \\
& = & g(\mathcal{G}(x)) \quad \mbox{ by definition of } g\\
& = & f(\mathcal{G}(x),\mathcal{G}(x)) \quad \mbox{ by representability of } g\\
& = & \mathcal{G}(\ulcorner\mathcal{G}(x)\urcorner) \quad \mbox{
by definition of } f.
\end{eqnarray*}$\Box$

\begin{apl}{G\"{o}del's First Incompleteness Theorem} Let $Prov(y,x)$
stand for ``$y$ is the G\"{o}del number of a proof of a statement
whose G\"{o}del number is $x$.'' Then let $$\mathcal{E}(x) \equiv
(\forall y)\neg Prov(y,x).$$ A fixed point for this
$\mathcal{E}(x)$ in a consistent and $\omega$-consistent theory is
a sentence that is equivalent to its own statement of
unprovability.
\end{apl}

\begin{apl}{G\"{o}del-Rosser's Incompleteness Theorem}
Let $Neg:\mathbb{N} \lra \mathbb{N}$ be defined for G\"{o}del
numbers as follows $$Neg(\ulcorner\mathcal{B}(x)\urcorner)\quad =
\quad \ulcorner\neg\mathcal{B}(x)\urcorner$$ Let $$\mathcal{E}(x)
\equiv (\forall y) (Prov(y,x)\rightarrow(\exists w) (w<y) \wedge
Prov(w, Neg(x))).$$ A fixed point for this $\mathcal{E}(x)$ in a
consistent theory is a sentence that is equivalent to its own
statement of unprovability.
\end{apl}

\begin{apl}{Tarski's Theorem} Let us assume that there exists a
well-formed formula $\mathcal{T}(x)$ that expresses the fact that
$x$ is the G\"{o}del number of a (true) theorem in the theory. Set
$$\mathcal{E}(x) \equiv \neg \mathcal{T}(x).$$ A fixed point of
$\mathcal{E}(x)$ shows that $\mathcal{T}(x)$ does not do what it
is supposed to do. We conclude that a theory in which the
diagonalization lemma holds cannot express its own theoremhood.
\end{apl}

\begin{apl}{Parikh Sentences} There are true sentences that have very long proofs, but
there are relatively short proof of the fact that the sentences
are provable. This amazing result about lengths of proofs can be
found on page 496 of R. Parikh's famous paper {\it Existence and
Feasibility in Arithmetic} \cite{MR46:3287}. Consider a consistent
theory that contains Peano Arithmetic. We shall deal with the
following predicates:
\begin{itemize}
\item $Prflen(m,x)\equiv$ $m$ is the length (in symbols) of a proof of
a statement whose G\"{o}del number is $x$.  This is decidable
because there are only a finite number of proofs of length $m$.
\item $P(x) \equiv \exists y Prov(y,x)$ i.e. there exists a proof
of a statement whose G\"{o}del number is $x$.
\item $\mathcal{E}_n(x) \equiv \neg(\exists m < n \quad
Prflen(m,x))$.
\end{itemize}
Applying the diagonalization lemma to $\mathcal{E}_n(x)$ gives us
a fixed point $\mathcal{C}_n$ such that
$$\vdash \mathcal{C}_n \longleftrightarrow
\mathcal{E}_n(\ulcorner\mathcal{C}_n\urcorner) \equiv \neg(\exists
m < n \quad Prflen(m,\ulcorner\mathcal{C}_n\urcorner)).$$ In other
words $\mathcal{C}_n$ says \newline\centerline{``I do not have a
proof of myself shorter than n.''}
 If $\mathcal{C}_n$ is false, then
there is a proof shorter than $n$ of $\mathcal{C}_n$ and the
system is not consistent.

Consider the following {\em short} proof of $P(\mathcal{C}_n)$
\begin{enumerate}
\item If $\mathcal{C}_n$  does not have any proof, then $\mathcal{C}_n$
is true.
\item If $\mathcal{C}_n$ is true, we can check all proofs of length
less than $n$ and prove $\mathcal{C}_n.$
\item From 1 and 2 we have that if $\mathcal{C}_n$  does not have
a proof, then we can prove $\mathcal{C}_n.$ i.e. $\neg
P(\mathcal{C}_n) \longrightarrow P(\mathcal{C}_n).$
\item $\therefore P(\mathcal{C}_n).$
\end{enumerate}

This proof can be formulated in Peano Arithmetic in a fairly short
proof. In contrast $n$ can be chosen to be fairly large. So we
have a statement $\mathcal{C}_n$ which has a very long proof, but
a short proof of the fact that it has a proof.
\end{apl}

\begin{apl}{L\"{o}b's Paradox} We prove that every logical
sentence is true. The standard notation for the G\"{o}del number
of a wff $\mathcal{C}$ is $\ulcorner \mathcal{C}\urcorner .$ In
contrast, if $n$ is an integer then we shall write $\llcorner n
\lrcorner$ for the wff that corresponds to that number. Obviously
$\llcorner\ulcorner \mathcal{C} \urcorner\lrcorner = \mathcal{C}$

Let $\mathcal{A}$ be any sentence. We shall prove that it is
always true. Use the diagonalization lemma on
$$\mathcal{E}(x) \equiv \llcorner x\lrcorner \Rightarrow \mathcal{A}.$$
A fixed point for this $\mathcal{E}(x)$ is a $\mathcal{C}$ such
that
$$\vdash \mathcal{C} \longleftrightarrow
\mathcal{E}(\ulcorner\mathcal{C}\urcorner) \equiv
(\llcorner\ulcorner\mathcal{C}\urcorner\lrcorner \Rightarrow
\mathcal{A}) = (\mathcal{C} \Rightarrow \mathcal{A}).$$ So
$\mathcal{C}$ is equivalent to $\mathcal{C} \Rightarrow
\mathcal{A}.$ Assume, for a second that $\mathcal{C}$ is true.
Then $\mathcal{C} \Rightarrow \mathcal{A}$ is also true. By modus
ponens $\mathcal{A}$ is also true. So by assuming $\mathcal{C}$ we
have proven $\mathcal{A}$. This is exactly what $\mathcal{C}
\Rightarrow \mathcal{A}$ says and hence it is true as is its
equivalent $\mathcal{C}$ and so $\mathcal{A}$ is true.

This looks like a real paradox. It seems to me that the paradox
arises because we did not put a restriction on the wffs
$\mathcal{E}(x)$ for which we are permitted to use the
diagonalization lemma. The L\"{o}b's paradox is related to Curry's
paradox which shows that we must restrict the comprehension scheme
in axiomatic set theory.
\end{apl}

Let us move from logic to computability theory. We shall use the
language and notation of \cite{MR81i:03001}.
\begin{teo}[The Recursion Theorem] Let $h:\mathbb{N}\lra
\mathbb{N}$ be a total computable function. There exists an $n_0
\in \mathbb{N}$ such that $$\phi_{h(n_0)} \quad = \quad
\phi_{n_0}.$$
\end{teo}
{\bf Proof.} Let $\mathcal{F}$ be the set of unary computable
functions. Consider $f:\mathbb{N}\times \mathbb{N}\longrightarrow
\mathcal{F}$ be defined as $f(m,n) \cong \phi_{\phi_n(m)}$. If
$\phi_n(m)$ is undefined, then $f(m,n)$ is also undefined. Letting
the operator $\Phi_h:\mathcal{F} \longrightarrow \mathcal{F}$ be
defined as $\Phi_h(\phi_n) = \phi_{h(n)}$. We have the following
square:
\sqr{\mathbb{N}}{\mathcal{F}}{f}{{\Phi_h}}{g}{\bigtriangleup} $g$
is defined as $g(m) = \phi_{h(\phi_m(m))}$. By the S-M-N theorem
there is a total computable function $s(m)$ such the
$\phi_{h(\phi_m(m))} = \phi_{s(m)}$. Since $s$ is total and
computable, there exists a number $t$ such that $s(m)=\phi_t(m)$
and so $g$ is representable because $g(m) =
\phi_{h(\phi_m(m))}=\phi_{s(m)}= \phi_{\phi_t(m)} = f(m,t)$. So
there is a fixed point of $\Phi_h$ at $n_0 = \phi_{\phi_t(t)}$.
Explicitly we have
\begin{eqnarray*}
\phi_{h({\phi_t(t)})} & = & \Phi_h(\phi_{\phi_t(t)})
\quad \mbox{ by definition of } \Phi_h \\
& = & \Phi_h(f(t,t)) \quad \mbox{ by definition of } f \\
& = & g(t) \quad \mbox{ by definition of } g\\
& = & f(t,t) \quad \mbox{ by representability of } g\\
& = & \phi_{\phi_t(t)} \quad \mbox{ by definition of } f.
\end{eqnarray*}$\Box$

\begin{apl}{Rice's Theorem} Every nontrivial property of computable functions is not decidable.
Let $\mathcal{A}$ be a nonempty proper subset of $\mathcal{F},$
the set of all unary computable functions. Let $A=\{x | \phi_x
\in \mathcal{A} \}.$ Then A is not recursive. We prove this by
assuming (wrongly) that $A$ is recursive. Let $a \in A$ and $b
\not\in A$. Define the function $h$ as follows.

$$ h(x)= \left \{
\begin{array}{r@{\quad : \quad}l}
a & \mbox{if } x \not\in A\\
b & \mbox{if } x \in A. \end{array}\right.$$

By definition $x \in A$ iff $h(x) \not\in A$. From our
assumption, we have that $h$ is computable (and total). Hence by
the recursion theorem, there is an $n_0$ such that $\phi_{h(n_0)}
= \phi_{n_0}$ Now we have the following contradiction:
\begin{eqnarray*}
n_0 \in A & \Longleftrightarrow & h(n_0) \not\in A
\quad \mbox{ by definition of } h \\
& \Longleftrightarrow & \phi_{h(n_0)}\not\in \mathcal{A} \quad \mbox{ by the definition of } A \\
& \Longleftrightarrow & \phi_{n_0} \not\in \mathcal{A} \quad \mbox{ by the recursion theorem } \\
& \Longleftrightarrow & n_0 \not\in A \quad \mbox{ by definition
of } A.
\end{eqnarray*}
\end{apl}

\begin{apl}{Von Neumann's Self-reproducing Machines} A
self-reproducing machine is a computable function that always
outputs its own description. It might seem impossible to construct
such a self-reproducing machine since in order to construct such
a machine, we would need to know its description and hence know
the machine in advance. However, by a simple application of the
recursion theorem, we get such a machine.

By a description of a machine, we could mean the number of the
computable function i.e. a self-reproducing machine is a function
$\phi_n(x) = n.$ for all input $x$.

Let $f:\mathbb{N} \times \mathbb{N}\lra \mathbb{N}$ be the
computable projection function $f(y,x) = y$. By the S-M-N theorem
there exists a total computable function $s$ such that
$\phi_{s(y)}(x) = f(y,x)= y$. From the recursion theorem, there
exists an $n$ such that $\phi_n(x) = \phi_{s(n)}(x) = f(n,x) = n$.
\end{apl}

\section{Future Directions}
There are many possible ways that we can go on with this work. We
shall list a few.

The general Cantor's theorem can be generalized further so that
even more phenomena can be encompassed by this one theorem. For
example what if we have two sets $Y$ and $Y'$ and there is a onto
function from $Y$ to $Y'$. What does this say about the
relationship between $f: T \times T \lra Y$ and $f': T \times T
\lra Y'$? We should get the concept of a paradox ``reduction''
from one paradox to another.

Rather than simply talking about sets and functions, perhaps we
should be talking about partial orders and order preserving maps.
With this generalization, we might be able to not only get fixed
point theorems but also {\it least} fixed point theorems. There
are many simple least fixed point theorems such as ones for
continuous maps of cpo's and Scott domains; Kripke's definition of
truth \cite{MR87f:03009} and the Knaster-Tarski theorem.

Some more thought must go into Richards and L\"{o}b's paradoxes.
Although we have stated their limitations, the paradoxes remain.
Perhaps we are not formulating them correctly or perhaps there is
something intrinsically problematic about these paradoxes.

There are many fixed point theorems throughout logic and
mathematics that are not of the type described in Sections 3 and
4. Can we in some sense characterize those fixed point theorems
that are self-referential?

It seems that the key component of the diagonalization lemma is
the existence of a recursive $D:\mathbb{N} \lra \mathbb{N}$ that
is defined for all $\mathcal{B}(x)$ as
$$ D(\ulcorner\mathcal{B}(x)\urcorner) =
\ulcorner\mathcal{B}(\ulcorner\mathcal{B}(x)\urcorner)\urcorner.$$
Similarly, in order to have the recursion theorem we needed the
S-M-N theorem. These two properties of systems are the key to the
fact that the systems can talk about themselves. Are these two
properties related to each other? More importantly, can we find
other key properties in systems that make self-reference possible?

In the introduction of this paper we talked of the lack of an onto
function $T \lra Y^T$ and we said that $Y$ may be thought of as
truth-values or properties of objects in $T$. Can we find a better
word for $Y$?  In Section 5 where we talked about an onto function
$Lind^1 \lra {Lind^0}^{Lind^1}$ where $Lind^i$ is the Lindenbaum
classes of formula with $i$ variables. In what sense is $Lind^0$
the truth-values or properties of $Lind^1$? We then went on to
talk about an onto function $\mathbb{N} \lra
\mathcal{F}^\mathbb{N}$ where $\mathcal{F}$ is the set of unary
computable functions. We used this onto function to prove The
Recursion Theorem. In what sense is $\mathcal{F}$ the truth-values
or the properties of $\mathbb{N}$?

As for more instances of our theorems, the field is wide open.
There are many paradoxical phenomena and fixed point theorems that
we have not talked about. Some of them might not be amenable to
our scheme and some might not be.

\begin{itemize}
\item There are many of the semantic paradoxes that we did not discuss.
The Berry paradox asks one to consider the sentence ``Let x be the
first number that can not be described by any sentence with less
than 200 characters.'' We just described such a number.
\item The Crocodile's Dilemma is an ancient paradox that is
a deviously cute self-referential paradox. A crocodile steals a
child and the mother of the child begs for the return of her
beloved baby. The crocodile responds "I will return the child if
and only if you correctly guess whether or not I will return your
child." The mother cleverly responds that he will keep the child.
What is an honest crocodile to do?!?

\item There is a belief that all paradoxes would melt away if there were no
self-referential statements. Yablo's Non-self-referential Liar's
Paradox was formulated counteract that thesis. There is a sequence
of statements such that none of them ever refer to themselves and
yet they are all both true and false. Consider the sequence
$$(S_i): \mbox{ For all } k>i, S_k \mbox{ is untrue.} $$
Suppose $S_n$ is true for some $n$. Then $S_{n+1}$ is false as are
all subsequent statements. Since all subsequent statements are
false, $S_{n+1}$ is true which is a contradiction. So in contrast,
$S_n$ is false for all $n$. That means that $S_1$ is true and
$S_2$ is true etc etc. Again we have a contradiction.
\item
Brandenburger's Epistemic Paradox \cite{brandenburger} considers
the situation where
\begin{quote}
{\it Ann believes that Bob believes that Ann believes that Bob has
a false belief about Ann.}
\end{quote}
Now ask yourself the following question: Does Ann believe that Bob
has a false belief about Ann?  With much thought, you can see that
this is a paradoxical situation.
\item The Ackermann function is not
a primitive recursive function. One hears the phrase that
Ackermann's function ``diagonalizes-out'' of primitive recursive
functions.
\item There is a famous Paris-Harrington result which says that
certain generalized Ramsey theorems can not be proven in Peano
Arithmatic. Kanamori and McAloon \cite{MR88i:03095} make the
connection to the Ackermann function. Just as the Ackermann
function ``diagonalized-out'' of primitive recursiveness, so too,
generalized Ramsey theory is ``diagonalized-out'' of Peano
Arithmetic. Both of these are really stating limitations of the
systems.
\end{itemize}

There are many instances of fixed point theorems that might be put
into the form of our scheme.
\begin{itemize}
\item Borodin's Gap Theorem is a type of fixed point theorem in
complexity theory that might be right for our scheme.
\item We again mention the Knaster-Tarski theorem about monotonic
functions between preorders. There is also a much used theorem
about fixed points of continuous functions between cpo's.
\item As the ultimate in self-reference, we would like to mention
Kripke's theory of truth that he used to banish self-referential
paradoxes. It is, in essence, a type of fixed point theorem. It
would really be nice to formulate that way of dealing with
paradoxes in our language.
\item Brouwer's fixed point theorem, or the far simpler
intermediate value theorem.
\item Nash's equilibria
theorem and its many generalizations from game theory.
\end{itemize}

There are several theorems from ``real'' mathematics that are
proved via diagonalization proofs. We might be able to put them
into our language.
\begin{itemize}
\item Baire's category theory about metric spaces.
\item Montel's theorem from complex function theory.
\item Ascoli theorem from topology.
\item Helly's theorem about limits of distributions.
\end{itemize}

The following ideas are a little more ``spacey.''
\begin{itemize}
\item G\"{o}del's second incompleteness theorem about the
unprovability within arithmetic of the consistency of arithmetic.
This theorem is a simple consequence of the first incompleteness
theorem. However Kreisal has a direct model theoretic proofs that
uses a diagonal method (see, e.g., page 860 of Smory\'nski's
article in \cite{handbook}.) This proof seems amenable to our
scheme.
\item Many of Chaitin's algorithmic information theory arguments
seem to fit our scheme.
\item We worked out G\"{o}del's first incompleteness theorem which
showed that (using the language of the introduction) arithmetic
can not completely talk about its own provability. What about
G\"{o}del's completeness theorem? Certain weak systems can
completely talk about their own provability. Can this be stated as
some type of fixed point theorem?
\end{itemize}

\bibliography{parbib.bib}

\begin{thebibliography}{10}

\bibitem{handbook}
{\em Handbook of mathematical logic}.
\newblock North-Holland Publishing Co., Amsterdam, 1977.
\newblock Edited by Jon Barwise, With the cooperation of H. J. Keisler, K.
  Kunen, Y. N. Moschovakis and A. S. Troelstra, Studies in Logic and the
  Foundations of Mathematics, Vol. 90.

\bibitem{MR52:16108}
Theodore Baker, John Gill, and Robert Solovay.
\newblock Relativizations of the {${\cal P}=?{\cal N}{\cal P}$} question.
\newblock {\em SIAM J. Comput.}, 4(4):431--442, 1975.

\bibitem{brandenburger}
Adam Brandenburger.
\newblock The power of paradox.
\newblock {\em available at
  http://www.people.hbs.edu/abrandenburger/paradox-03-12-021.pdf}.

\bibitem{MR81i:03001}
Nigel Cutland.
\newblock {\em Computability}.
\newblock Cambridge University Press, Cambridge, 1980.
\newblock An introduction to recursive function theory.

\bibitem{MR87f:03009}
Melvin Fitting.
\newblock Notes on the mathematical aspects of {K}ripke's theory of truth.
\newblock {\em Notre Dame J. Formal Logic}, 27(1):75--88, 1986.

\bibitem{Heller}
Alex Heller.
\newblock An existence theorem for recursion categories.
\newblock {\em J. Symbolic Logic}, 55(3):1252--1268, 1990.

\bibitem{MR83j:68002}
John~E. Hopcroft and Jeffrey~D. Ullman.
\newblock {\em Introduction to automata theory, languages, and computation}.
\newblock Addison-Wesley Publishing Co., Reading, Mass., 1979.
\newblock Addison-Wesley Series in Computer Science.

\bibitem{MR91g:18007}
Hagen Huwig and Axel Poign{\'e}.
\newblock A note on inconsistencies caused by fixpoints in a {C}artesian closed
  category.
\newblock {\em Theoret. Comput. Sci.}, 73(1):101--112, 1990.

\bibitem{MR88i:03095}
Akihiro Kanamori and Kenneth McAloon.
\newblock On {G}\"odel incompleteness and finite combinatorics.
\newblock {\em Ann. Pure Appl. Logic}, 33(1):23--41, 1987.

\bibitem{MR95k:00009}
Shaughan Lavine.
\newblock {\em Understanding the infinite}.
\newblock Harvard University Press, Cambridge, MA, 1994.

\bibitem{MR39:4075}
F.~William Lawvere.
\newblock Diagonal arguments and cartesian closed categories.
\newblock In {\em Category Theory, Homology Theory and their Applications, II
  (Battelle Institute Conference, Seattle, Wash., 1968, Vol. Two)}, pages
  134--145. Springer, Berlin, 1969.

\bibitem{setsformath}
F.~William Lawvere and Robert Rosebrugh.
\newblock {\em Sets for Mathematics}.
\newblock Cambridge University Press.

\bibitem{MR93m:18001}
F.~William Lawvere and Stephen~H. Schanuel.
\newblock {\em Conceptual mathematics}.
\newblock Buffalo Workshop Press, Buffalo, NY, 1991.
\newblock A first introduction to categories, With the assistance of Emilio
  Faro, Fatima Fenaroli and Danilo Lawvere.

\bibitem{MR2001j:18001}
Saunders Mac~Lane.
\newblock {\em Categories for the working mathematician}, volume~5 of {\em
  Graduate Texts in Mathematics}.
\newblock Springer-Verlag, New York, second edition, 1998.

\bibitem{Manin}
Yuri~I. Manin.
\newblock Classical computing, quantum computing, and {S}hor's factoring
  algorithm.
\newblock {\em Ast\'erisque}, (266):Exp.\ No.\ 862, 5, 375--404, 2000.
\newblock S\'eminaire Bourbaki, Vol. 1998/99.

\bibitem{MR99b:03002}
Elliott Mendelson.
\newblock {\em Introduction to mathematical logic}.
\newblock Chapman \& Hall, London, fourth edition, 1997.

\bibitem{MR91c:03049}
Philip~S. Mulry.
\newblock Categorical fixed point semantics.
\newblock {\em Theoret. Comput. Sci.}, 70(1):85--97, 1990.
\newblock Fourth Workshop on Mathematical Foundations of Programming Semantics
  (Boulder, CO, 1988).

\bibitem{MR46:3287}
Rohit Parikh.
\newblock Existence and feasibility in arithmetic.
\newblock {\em J. Symbolic Logic}, 36:494--508, 1971.

\bibitem{MR93i:03021}
Du{\v{s}}ko Pavlovi{\'c}.
\newblock On the structure of paradoxes.
\newblock {\em Arch. Math. Logic}, 31(6):397--406, 1992.

\bibitem{MR91k:03018}
Andrew~M. Pitts and Paul Taylor.
\newblock A note on {R}ussell's paradox in locally {C}artesian closed
  categories.
\newblock {\em Studia Logica}, 48(3):377--387, 1989.

\end{thebibliography}
\bibliographystyle{plain}
\nocite{*}

\noindent Department of Computer and Information Science\\
Brooklyn College, CUNY\\
Brooklyn, N.Y. 11210\\
\\
and\\
\\
Department of Computer Science\\
The Graduate Center, CUNY \\
365 Fifth Avenue\\
New York, N.Y. 10016\\
\\
email:noson@sci.brooklyn.cuny.edu
\end{document}